\documentclass[a4paper,leqno]{amsart}
\usepackage{graphicx, amsmath, amssymb, amsthm, hyperref, verbatim, multicol, mathrsfs, tikz, caption, subcaption}

\theoremstyle{definition}
\newtheorem{thm}{Theorem}[section]

\newtheorem{lemm}[thm]{Lemma}
\newtheorem{cor}[thm]{Corollary}

\newtheorem{prop}[thm]{Proposition}

\DeclareMathOperator{\Lip}{Lip}

\DeclareMathOperator{\diam}{diam}

\allowdisplaybreaks

\begin{document}

\title{Singular quasisymmetric mappings in dimensions two and greater}
\author{Matthew Romney}
\address{Department of Mathematics and Statistics, University of Jyv\"askyl\"a, P.O. Box 35 (MaD), FI-40014, University of Jyv\"askyl\"a, Finland.}
\email{matthew.d.romney@jyu.fi}
\keywords{Quasiconformal mapping, metric space, absolute continuity}

\thanks{This research was supported by the Academy of Finland grant 288501 and by the ERC Starting Grant 713998 GeoMeG. \newline {\it 2010 Mathematics Subject Classification.} 30L10.}

\begin{abstract}
For all $n \geq 2$, we construct a metric space $(X,d)$ and a quasisymmetric mapping $f\colon [0,1]^n \rightarrow X$ with the property that $f^{-1}$ is not absolutely continuous with respect to the Hausdorff $n$-measure on $X$. That is, there exists a Borel set $E \subset [0,1]^n$ with Lebesgue measure $|E|>0$ such that $f(E)$ has Hausdorff $n$-measure zero. The construction may be carried out so that $X$ has finite Hausdorff $n$-measure and $|E|$ is arbitrarily close to 1, or so that $|E| = 1$. This gives a negative answer to a question of Heinonen and Semmes.
\end{abstract}

\maketitle

\section{Introduction} \label{sec:introduction}

Let $(X,d_X,\mu)$ and $(Y, d_Y, \nu)$ be metric measure spaces. For our purposes, a {\it metric measure space} is a metric space equipped with a Borel measure which assigns a positive (though possibly infinite) value to all metric balls.  A (homeomorphic) mapping $f\colon X \rightarrow Y$ is {\it absolutely continuous in measure} if for all measurable sets $E \subset X$, $\mu(E) = 0$ implies $\nu(f(E)) = 0$. A mapping $f\colon X \rightarrow Y$ is {\it quasisymmetric} if there exists a homeomorphism $\eta\colon [0, \infty) \rightarrow [0, \infty)$ such that
$$\frac{d_Y(f(x),f(y))}{d_Y(f(x),f(z))} \leq \eta\left( \frac{d_X(x,y)}{d_X(x,z)} \right)$$
for all triples of distinct points $x,y,z \in X$. Notice that this definition only depends on the metrics $d_X$ and $d_Y$. The class of quasisymmetric mappings form a natural and well-studied class of geometry-preserving mappings, with applications in fields such as complex dynamics and geometric group theory. See \cite{TV} and \cite[Ch. 10-11]{hei:lectures} for the basic theory.

A fundamental problem is to determine which metric spaces $(X,d_X)$ are quasisymmetrically equivalent to some Euclidean space $\mathbb{R}^n$ (or the sphere $\mathbb{S}^n$). This is the {\it quasisymmetric uniformization problem} for Euclidean spaces (alternatively, one can speak of characterizing the {\it conformal gauge} of $\mathbb{R}^n$ \cite[Ch. 15]{hei:lectures}). In full generality, this is a difficult problem. Often one can make progress by imposing some kind of geometric condition on the space $X$. For $n=2$, a simple and elegant geometric characterization was obtained by Bonk and Kleiner \cite{BonkKleiner} under the assumption that $X$ is Ahlfors 2-regular. The Bonk--Kleiner theorem states that an Ahlfors 2-regular metric sphere $X$ is quasisymmetrically equivalent to the standard sphere if and only if it is linearly locally contractible (that is, there exists $\lambda\geq 1$ such that every ball of radius $r<\diam(X)/\lambda$ contracts to a point inside a ball of radius $\lambda r$). 
We recall that $X$ is {\it Ahlfors $n$-regular} if there exists a $C\geq 1$ such that $C^{-1}r^n \leq \mathcal{H}^n(B(x,r)) \leq Cr^n$ for all metric balls $B(x,r) \subset X$ with $r < \diam X$. Here, $\mathcal{H}^n$ denotes the $n$-dimensional Hausdorff measure on $X$. For some other contributions to the quasisymmetric uniformization problem, see \cite{Laa:02}, \cite{LoRajRa:18}, \cite{Raj16}, \cite{Sem:96b}, \cite{Wu:grushin}.

Without a strong geometric condition like Ahlfors $n$-regularity, the situation is far less nice. A long-standing open question in this area is whether a quasisymmetric mapping $f\colon (X,d) \rightarrow \mathbb{R}^n$, where $n \geq 2$ and $(X,d)$ is some metric space homeomorphic to $\mathbb{R}^n$, must be absolutely continuous with respect to Hausdorff $n$-measure on $X$ (denoted throughout this paper by $\mathcal{H}^n$). This question appears most prominently as Questions 15 and 16 in the survey article of Heinonen and Semmes \cite{HeiS:97}. Since $\mathbb{R}^n$ is often thought of as the source space, the problem is frequently called the {\it inverse absolute continuity problem} for quasisymmetric mappings. The property of absolute continuity in measure is also referred to as {\it Lusin's condition (N)} in the literature. 

The purpose of this paper is to give a negative answer to this problem, as stated in Questions 15 and 16 in \cite{HeiS:97}. Before giving a precise statement of our results, we offer a few general remarks.

The Heinonen--Semmes problem is part of a broader set of questions related to the regularity of quasiconformal mappings and their generalizations. The classical problem originated from work of Gehring in the 1970s \cite{Gehr:75}, \cite{Gehr:76}. It asks whether, given an $n$-dimensional hyperplane (or hypersurface) $V$ in $\mathbb{R}^{n+1}$ ($n\geq 2$), a quasiconformal homeomorphism of $\mathbb{R}^{n+1}$ can map a subset $E \subset V$ of positive $n$-measure onto a set of Hausdorff $n$-measure zero. In \cite[Qu. 5.10]{Vai:81}, V\"ais\"al\"a asks the analogous question for the case of a quasisymmetric embedding $f\colon \mathbb{R}^n \rightarrow \mathbb{R}^N$, $N >n$. These problems are still open, although we hope that the results of the present paper can inspire an approach to their solution. [Added in October 2018: the author and D. Ntalampekos, in the preprint \cite{NR:18}, have resolved the problems of Gehring and V\"ais\"al\"a for the case $n=2$.]

One may ask about absolute continuity in the other direction. A result of Tyson \cite[Cor. 5.10]{Tys:00} states that a locally quasisymmetric mapping $f\colon \mathbb{R}^n \rightarrow X$ ($n \geq 2$) is absolute continuous in measure (as a mapping from $\mathbb{R}^n$ to $X$), provided that $X$ has locally finite Hausdorff $n$-measure. Tyson also shows that the Hausdorff $n$-measure on $X$ cannot be too small, for any metric space $X$ quasisymmetrically equivalent to $\mathbb{R}^n$. Precisely, there exists a constant $C>0$ depending on $n$ and the quasisymmetry control function $\eta$ such that $\mathcal{H}^n(B(x,r)) \geq Cr^n$ for all $x \in X$ and sufficiently small $r>0$ \cite[Cor. 3.10]{Tys:00}. These results generalize earlier work of V\"ais\"al\"a in \cite{Vai:81} for the case that $X$ is a subset of Euclidean space. 

A similar question is whether the boundary extension of a quasiconformal mapping from the unit ball $\mathbb{B}^n$ in $\mathbb{R}^n$ ($n \geq 3$) to a Jordan domain in $\mathbb{R}^n$ with $(n-1)$-rectifiable boundary must be absolutely continuous with respect to $(n-1)$-dimensional Hausdorff measure. Examples where boundary absolute continuity fails were given by Heinonen in \cite{Hei:96}. The topic of boundary extensions of quasiconformal mappings was further studied by Astala, Bonk, and Heinonen in \cite{ABH:02}, partly motivated by the inverse absolute continuity problem for quasiconformal mappings on hyperplanes. For further discussion on these and related problems, see \cite[Sec. 6]{Hei:96}, \cite[Sec. 1]{ABH:02}, \cite[Sec. 17]{Raj16}. 


We continue now with the statement of the main theorems. Since the problem at hand is essentially local in nature, there is no loss of generality by considering maps defined on the cube $[0,1]^n$ instead of the whole space $\mathbb{R}^n$. We give two versions of our result. In the first, there is no requirement that the Hausdorff $n$-measure be locally finite. In this case, inverse absolute continuity can fail especially badly. 

\begin{thm} \label{thm:main} 
Let $n \geq 2$ and $\alpha \in (0,n)$. There exists a metric space $(X,d)$ and a quasisymmetric homeomorphism $f\colon [0,1]^n \rightarrow X$ with the property that $f$ maps a Borel set $E \subset [0,1]^n $ of full measure onto a subset of $X$ with Hausdorff dimension at most $\alpha$. 
\end{thm}  

Observe that the mapping $f$ in Theorem \ref{thm:main} takes the null set $[0,1]^n \setminus E$ onto a set of full Hausdorff $n$-measure in the target. In particular, we see from this that any space $X$ satisfying the conclusions of Theorem \ref{thm:main} cannot have locally finite Hausdorff $n$-measure. If so, this would contradict the above-mentioned theorem of Tyson on the absolute continuity of quasisymmetric mappings from $\mathbb{R}^n$ for $n \geq 2$. 



The second version of our theorem does yield a metric space with locally finite Hausdorff measure. 



\begin{thm} \label{thm:main2} 
	Let $n \geq 2$, $\alpha \in (0,n)$ and $\delta \in (0,1)$. There exists a metric space $(X,d)$ of finite Hausdorff $n$-measure and a quasisymmetric homeomorphism $f\colon [0,1]^n \rightarrow X$ with the property that $f$ maps a Borel set $E$ of Lebesgue measure $|E| \geq \delta$ onto a subset of $X$ with Hausdorff dimension at most $\alpha$. 
\end{thm}  

The majority of the paper will be devoted to a construction of the space $(X,d)$ described in Theorem \ref{thm:main} and proving the required properties. Theorem \ref{thm:main2} is proved by the same construction, with a certain modification to keep the Hausdorff measure of $X$ finite.


Our construction scheme applies to all Euclidean dimensions two and greater. For the one-dimensional version of the problem, it is a classical fact that one-dimensional quasisymmetric mappings may be highly singular, even for quasisymmetric homeomorphisms of the real line. This was originally discovered by Ahlfors and Beurling in their investigations of boundary behavior of quasiconformal mappings \cite{AhlBeu:56}. The extent to which a quasisymmetric mapping of the real line can distort Hausdorff dimension was the subject of a later paper of Tukia \cite{Tuk:89}. His main result states that for all $\alpha>0$, there exists a quasisymmetric mapping $f\colon [0,1] \rightarrow [0,1]$ with the property that there exists a subset $E \subset [0,1]$ such that both $f(E)$ and $[0,1] \setminus E$ have Hausdorff dimension at most $\alpha$. Indeed, the construction in this paper was partly inspired by Tukia's construction. In our theorems, we obtain the same type of control on the Hausdorff dimension of $f(E)$. However, we do not know to what extent the Hausdorff dimension of $[0,1]^n \setminus E$ in Theorem \ref{thm:main} can be reduced.  

On the other hand, the inverse absolute continuity property for a quasisymmetric mapping $f\colon \mathbb{R}^n \rightarrow X$ is known to hold in the $n \geq 2$ case if $X$ satisfies additional assumptions. For instance, it holds if $X$ is Ahlfors $n$-regular \cite{HeiS:97}. In the two-dimensional case, Rajala has recently shown the answer to also be affirmative if $X$ satisfies a reciprocal upper bound on the modulus of path families associated to a quadrilateral (Proposition 17.2 in \cite{Raj16}). 

Finally, we highlight a pair of applications of our main results. The first is to the related topic of quasiconformal uniformization. This gives a negative answer to Question 17.3 of \cite{Raj16}. 
\begin{cor} \label{cor:qc}
There exists a metric space $(X,d)$ of locally finite Hausdorff 2-measure and a quasisymmetric mapping $f\colon [0,1]^2 \rightarrow X$ which fails to be quasiconformal.
\end{cor}
By {\it quasiconformal}, we refer to the geometric definition based on the notion of modulus of a path family. Corollary \ref{cor:qc} follows from the fact that inverse absolute continuity is known to hold for quasiconformal mappings under the given hypotheses \cite[Rem. 8.3]{Raj16}. Interestingly, one can construct a metric space $X$ with locally finite Hausdorff 2-measure admitting a quasiconformal mapping $f\colon \mathbb{R}^2 \rightarrow X$ which fails to the absolutely continuous in measure (as a mapping from $\mathbb{R}^2$ to $X$) \cite[Prop. 17.1]{Raj16}. As already mentioned, such a mapping cannot be quasisymmetric. In summary, there is no general relationship between the existence of a quasiconformal parametrization and the existence of a quasisymmetric parametrization in this setting. 

The second application is to the notion of {\it conformal dimension of measures}, recently introduced by Bate and Orponen \cite{BaOr:18}. This is a measure-theoretic counterpart to the conformal dimension of a metric space (see the book of Mackay and Tyson \cite{MacTy:10}). The {\it conformal dimension} of a locally finite Borel measure $\mu$ on a metric space is the infimal value of the Hausdorff dimension of the image of a set of full $\mu$-measure under a quasisymmetric mapping onto any metric space. Expressed in this language, Theorem \ref{thm:main} reads as follows. 
\begin{cor}\label{cor:conformal_dimension}
The conformal dimension of Lebesgue measure on $\mathbb{R}^n$ is zero.
\end{cor}
As discussed in Question 1 of \cite{BaOr:18}, it is not clear whether there exists any measure with positive and finite conformal dimension. Corollary \ref{cor:conformal_dimension} supports the assertion that they do not exist. See also work of Meyer \cite{Mey:09}, where a similar concept is studied under the name {\it elliptic harmonic measure}.

The paper is organized as follows. In Section \ref{sec:construction}, we give a description of our main construction, a sequence of conformal weights $\rho_k$ on $[0,1]^n$ which yield a sequence of metric spaces $(X_k, d_k)$. These metrics converge pointwise to a metric $d$ on $[0,1]^n$. We take as $f$ the identity mapping. Section \ref{sec:estimates} contains basic estimates of the metrics $d_k$ at different scales and relative distances. In Section \ref{sec:quasisymmetry}, we verify that the mapping $f$ is indeed quasisymmetric. In Section \ref{sec:absolute_continuity}, we show that $f^{-1}$ is not absolutely continuous in measure. 
In Section \ref{sec:theorem2_proof}, we show how to modify our construction to obtain Theorem \ref{thm:main2}. 

\section{The construction} \label{sec:construction}

We take as our source space the unit $n$-cube $Q=[0,1]^n$ in $\mathbb{R}^n$. We use an iterative construction to define a sequence of metric spaces $(X_k, d_k)$ converging to a limit space $(X,d)$, which we show to be quasisymmetrically equivalent to $Q$. The metric on $X_k$ is specified via a conformal weight $\rho_k$ on $Q$.

Before beginning, we take some time to fix notation. Let $M \in \mathbb{N}$ and $L>1$ be sufficiently large. The admissible values of $M$ and $L$ depends on the dimension $n$ and the desired bound on the Hausdorff dimension of $f(E)$. This will be discussed in Section \ref{sec:absolute_continuity}. For the case $n=2$, we may take any values $M,L \geq 8$, provided we are not concerned with a particular bound on the Hausdorff dimension of $f(E)$. 

For each $k \in \mathbb{Z}_{\geq 0}$, let $\mathcal{V}_k$ denote the set of points $(i_1M^{-k},\ldots, i_nM^{-k})$ for $0 \leq i_1, \ldots, i_n \leq M^k$, and let $\mathcal{V} = \bigcup_k \mathcal{V}_k$. Write $I = (i_1, \ldots, i_n)$. For a given $k \in \mathbb{Z}_{\geq 0}$ and $1 \leq i_1, \ldots, i_n \leq M^k$, let $Q_k(I)$ denote the $n$-cube (or {\it cube}) of side length $M^{-k}$  whose vertices are $((i_1-\varepsilon_1)M^{-k}, \ldots, (i_k-\varepsilon_n)M^{-k})$, where $\varepsilon_1, \ldots, \varepsilon_n$ range over $\{0,1\}$, and let $q_k(I)$ denote the center point of this cube. Let $\mathcal{Q}_k$ denote the collection of all such cubes corresponding to the level $k$.  For a point $x \in Q$, let $Q_k(x)$ denote some cube $Q_k(I)$ containing $x$. Note that $x$ may lie on the boundary of multiple cubes, so $Q_k(x)$ is not uniquely determined. Unless otherwise noted, this choice can be arbitrary. Let $\mathcal{E}_k$ denote the set of points lying in $\partial Q$ for some $Q \in \mathcal{Q}_k$.

Next, we define the following collections relative to a given cube $Q_k(I)$. For its statement, we use the notation $d_{|\cdot|}(\cdot, \cdot)$ to denote Euclidean distance between sets. That is, $d_{|\cdot|}(E,F) = \inf\{|x-y|: x \in E, y \in F\}$. 

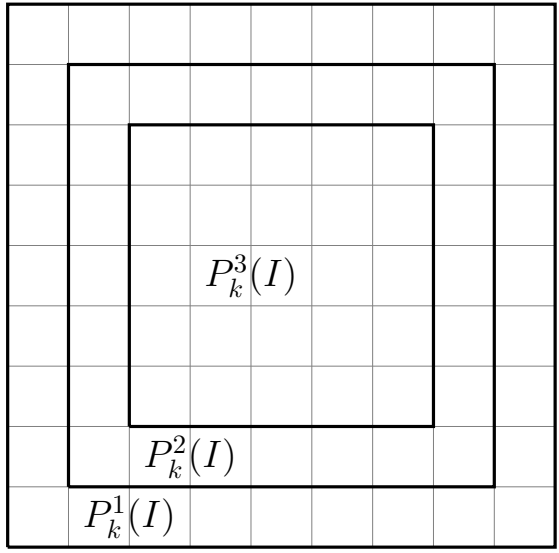
\begin{figure}[tb] \centering
		\begin{tikzpicture}[>=latex, scale=.8]
		\foreach \i in {1,...,8} {
			\draw [very thin,gray] (\i,0) -- (\i,9);
			\draw [very thin,gray] (0,\i) -- (9,\i);
		}
		\draw[very thick] (0,0) -- (9,0) -- (9,9) -- (0,9) -- (0,0);
		\draw[very thick] (1,1) -- (8,1) -- (8,8) -- (1,8) -- (1,1);
		\draw[very thick] (2,2) -- (7,2) -- (7,7) -- (2,7) -- (2,2); 
		\node at (2,.5) {\LARGE $P_k^1(I)$};
		\node at (3,1.5) {\LARGE $P_k^2(I)$};
		\node at (4,4.5) {\LARGE $P_k^3(I)$};
		\end{tikzpicture}
		\caption{The sets $P_k^1(I)$, $P_k^2(I)$, $P_k^3(I)$, $n=2$}
		\label{fig:construction}
\end{figure}

\begin{enumerate}
\item $Q_k^*(I) = \bigcup \{Q \in \mathcal{Q}_k: Q \cap Q_k(I) \neq \emptyset\}$.
\item $\mathcal{P}_k^1(I) = \{Q\in \mathcal{Q}_{k+1}: Q \subset Q_k(I) \text{ and } Q \cap \partial Q_k(I) \neq \emptyset\}$. 
\item $\mathcal{P}_k^2(I) = \{Q \in \mathcal{Q}_{k+1}: Q \subset Q_k(I), Q \notin \mathcal{P}_k^1(I), \text{ and } d_{|\cdot|}(Q, \bigcup \mathcal{P}_k^1(I)) < (n-1)M^{-(k+1)}\}$. 
\item $\mathcal{P}_k^3(I) = \{Q \in \mathcal{Q}_{k+1}: Q \subset Q_k(I), Q \notin \mathcal{P}_k^1(I) \cup \mathcal{P}_k^2(I) \}$. 
\end{enumerate}  
Observe that $|\mathcal{P}_k^1(I)| = M^n - (M-2)^n$, $|\mathcal{P}_k^2(I)| = (M-2)^n-(M-2n)^n$, and $|\mathcal{P}_k^3(I)| = (M-2n)^n$ for all cubes $Q_k(I)$. Next, let $P_k^1(I) = \bigcup \mathcal{P}_k^1(I)$, $P_k^3(I) = \bigcup \mathcal{P}_k^3(I)$ and $P_k^2(I) = \bigcup \mathcal{P}_k^2(I) \setminus (P_k^1(I) \cup P_k^3(I))$. See Figure \ref{fig:construction} for a representation of these sets in the 2-dimensional case. 

Define the weight $\rho_0\colon Q \rightarrow [0, \infty)$ by $\rho_0(x) = 1$ for all $x \in Q$. For the inductive step, given a weight $\rho_k$, define $\rho_{k+1}\colon Q \rightarrow [0, \infty)$ on the interior of each cube $Q_{k+1}(I)$ by the formula
\begin{equation}  \label{equ:rho_definition}
\rho_{k+1}(x) = \left\{ \begin{array}{cc} \rho_k(q_k(I)) & \text{ if } x \in P_k^1(I) \\ \rho_k(q_k(I))(M-2n+1) & \text{ if } x \in P_k^2(I) \\ \rho_k(q_k(I))/L & \text{ if } x \in P_k^3(I) \end{array} \right.  . 
\end{equation} 
Extend the definition to all points in $Q$ by requiring $\rho_{k+1}$ to be lower semicontinuous. 

There are three desired attributes which motivate the above definition of $\rho_k$. First, $\rho_{k+1}$ is much smaller than $\rho_k$ on the large set $P_k^3(I)$. This leads to the failure of inverse absolute continuity. Second, $\rho_{k+1}$ is much larger than $\rho_k$ on the set $P_k^2(I)$. The result of this is that the large-scale behavior of the resulting metric is not affected by later iterations of the weight, and in particular that the limiting metric is not degenerate. Third, $\rho_k$ is unchanged on the set $P_k^1(I)$. This ensures that the quasisymmetry condition is not upset by the iteration procedure. The precise values given in the $x \in P_k^2(I)$ and $x\in P_k^3(I)$ cases of \eqref{equ:rho_definition} are not essential, provided they fall within a certain range that achieves these three goals. 

In the following, let $r_k(I) = \rho_k(q_k(I))$ and $R = R(M,L) = L(M-2n+1)$. Given $x \in Q$, let $q_k(x)$ be the center point of $Q_k(x)$ and let $r_k(x) = \rho_k(q_k(x))$. We make the observation that for all cubes $Q_k(I)$ and all $Q_k(I') \subset Q_k^*(I)$,
\begin{equation} \label{equ:rho_comparability} 
\frac{1}{R} \leq \frac{r_k(I)}{r_k(I')} \leq R.
\end{equation}
 
The weights $\rho_k$ each yield a metric space $(X_k, d_k)$ in the usual fashion, by defining
$$d_k(x,y) = \inf \int_\gamma \rho_k \circ \gamma \,ds,$$
where the infimum is taken over all absolutely continuous paths in $Q$ with initial point $x$ and terminal point $y$ (for the remainder of the paper, the term {\it path} means absolutely continuous path). Let $\ell_k(\gamma)$ denote the $d_k$-length of the path $\gamma$. We also use $|\gamma|$ to denote the image of $\gamma$.

Informally, we define $d$ as the pointwise limit of the metrics $d_k$, and let $X$ be the resulting metric space. We postpone a precise definition until after a detailed analysis of the metrics $d_k$. 

\section{Basic estimates} \label{sec:estimates}

This section is dedicated to a number of basic lemmas describing the behavior of the metrics $d_k$. The first lemma states that passing from the weight $\rho_k$ to the weight $\rho_{k+1}$ does not decrease the length of a path which connects the boundary of a cube in $\mathcal{Q}_k$.

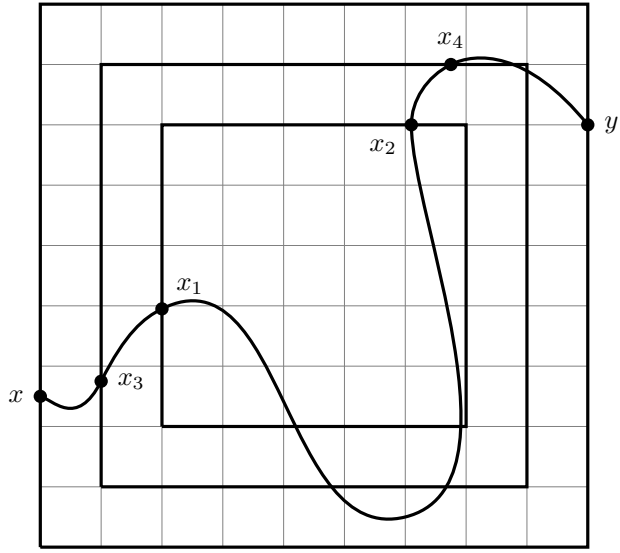
\begin{figure}[tb] \centering
		\begin{tikzpicture}[>=latex, scale=.8]
		\foreach \i in {1,...,8} {
			\draw [very thin,gray] (\i,0) -- (\i,9);
			\draw [very thin,gray] (0,\i) -- (9,\i);
		}
		\draw[very thick] (0,0) -- (9,0) -- (9,9) -- (0,9) -- (0,0);
		\draw[very thick] (1,1) -- (8,1) -- (8,8) -- (1,8) -- (1,1);
		\draw[very thick] (2,2) -- (7,2) -- (7,7) -- (2,7) -- (2,2); 
		\draw[very thick] (0,2.5) .. controls (.2,2.45) and (.6,1.95) .. (1, 2.75) .. controls (1.4,3.55) and (1.7,3.8) .. (2,3.95) .. controls (4, 4.95) and (4, 0) .. (6,.5) .. controls (8,1) and (6.1,5.5) .. (6.1,7) .. controls (6.1,7.4) and (6.35,7.8) .. (6.75, 8) .. controls (7.15,8.2) and (8,8.25) .. (9,7);
		\node[circle,fill=black,inner sep=0pt,minimum size=5pt,label=left:{$x$}] at (0, 2.5) {};
		\node[circle,fill=black,inner sep=0pt,minimum size=5pt,label=right:{$x_3$}] at (1, 2.75) {};
		\node[circle,fill=black,inner sep=0pt,minimum size=5pt,label=above right:{$x_1$}] at (2,3.95) {};
		\node[circle,fill=black,inner sep=0pt,minimum size=5pt,label=below left:{$x_2$}] at (6.1,7) {};
		\node[circle,fill=black,inner sep=0pt,minimum size=5pt,label=above:{$x_4$}] at (6.75, 8) {};
		\node[circle,fill=black,inner sep=0pt,minimum size=5pt,label=right:{$y$}] at (9, 7) {};
		\end{tikzpicture}
		\caption{The path $\gamma$ in Lemma \ref{lemm:length_of_curves}}
		\label{fig:path_gamma}
\end{figure}

\begin{lemm} \label{lemm:length_of_curves} 
Let $\gamma\colon [0,1] \rightarrow Q_k(I)$ be a path with endpoints $x,y \in \partial Q_k(I)$, for some cube $Q_k(I) \in \mathcal{Q}_k$. Then $\ell_{k+1}(\gamma) \geq d_k(x,y)$. 
\end{lemm}
\begin{proof}
If $\gamma$ is contained in $Q_k(I) \setminus P_k^3(I)$, then it is clear that $\ell_{k+1}(\gamma) \geq \ell_k(\gamma) \geq d_k(x,y)$. For the case that $|\gamma| \cap P_k^3(I) \neq \emptyset$, let $t_1$ denote the smallest value such that $\gamma(t_1) \in P_k^3(I)$ and let $t_2$ denote the largest value such that $\gamma(t_2) \in P_k^3(I)$. Similarly, let $t_3$ denote the largest value less than $t_1$ such that $\gamma(t_3) \in P_k^1(I)$ and let $t_4$ denote the smallest value greater than $t_2$ such that $\gamma(t_4) \in P_k^1(I)$. For $j \in \{1,\ldots, 4\}$, let $x_j = \gamma(t_j)$, and let $\widetilde{\gamma} = \gamma|[t_3, t_4]$. See Figure \ref{fig:path_gamma} for an illustration of a typical path $\gamma$. Then 
\begin{align} \label{equ:sublength}
\ell_{k+1}(\widetilde{\gamma}) & \geq r_k(I)\left(\frac{|x_1 - x_2|}{L} + (M-2n+1)(|x_1 - x_3| + |x_2 - x_4|) \right).
\end{align}  
We claim that the right-hand side of \eqref{equ:sublength} is bounded below by 
\begin{equation} \label{equ:sublength_difference} 
r_k(I)(|x_1 - x_2| + |x_1 - x_3| + |x_2 - x_4|).
\end{equation}  
Subtracting \eqref{equ:sublength_difference} from the right-hand side of \eqref{equ:sublength} and simplifying, we want to verify the inequality
\begin{equation*} 
(M-2n)(|x_1 - x_3| + |x_2 - x_4|) - \left(1 - \frac{1}{L}\right)|x_1 - x_2| \geq 0 .
\end{equation*} 
Since $|x_1 - x_3| \geq (n-1)M^{-1-k}$, $|x_2 - x_4| \geq (n-1)M^{-1-k}$ and $|x_1 - x_2| \leq n(M-2n)M^{-1-k}$, this reduces to 
$$2(n-1)(M-2n) - n(M-2n)\left(1- \frac{1}{L}\right) \geq 0,$$ 
which is satisfied. We conclude that $\ell_{k+1}(\widetilde{\gamma}) \geq r_k(I)(|x_1 - x_2| + |x_1 - x_3| + |x_2 - x_4|)$. This implies that $\ell_{k+1}(\gamma) \geq \ell_k(\gamma')$, where $\gamma'$ is the polygonal path from $x$ to $x_3$ to $x_1$ to $x_2$ to $x_4$ to $y$. Since $\ell_k(\gamma') \geq d_k(x,y)$, the lemma now follows. 
\end{proof}

As a consequence of Lemma \ref{lemm:length_of_curves}, the sequence $(d_m(x,y))$ is eventually nondecreasing for all $x,y \in \bigcup_k \mathcal{E}_k$. This is the content of the next lemma.

\begin{lemm} \label{lemm:increasing_property}
Suppose that $x, y \in \mathcal{E}_k$ for some $k \in \mathbb{Z}_{\geq 0}$. Then for all $m \geq k$, $d_m(x,y) \geq d_k(x,y)$.
\end{lemm}
\begin{proof} 
We begin with the case that $m = k+1$. Let $\gamma\colon [0,1] \rightarrow Q$ be a path from $x$ to $y$. Write $[0,1]$ as the countable union of intervals $[a_j,b_j]$ which do not overlap except at the endpoints, such that the image of $\gamma|[a_j,b_j]$ is contained in some cube $Q_k(I_j)$. Then $\ell_{k+1}(\gamma) = \sum_j \ell_{k+1}(\gamma|[a_j,b_j]) \geq \sum_j d_k(\gamma(a_j), \gamma(b_j)) \geq d_k(x,y)$. Infimizing over all possible paths shows that $d_{k+1}(x,y) \geq d_k(x,y)$. 

Since $x,y \in \mathcal{E}_{k+1}$, the result now follows by induction. 
\end{proof} 

The next lemma shows that the metrics $d_m$ are uniformly bounded on a cube $Q_k(I) \in \mathcal{Q}_k$ for all $m \geq k$. It also implies that the metrics $d_m$ are equicontinuous at a point $x \in \mathcal{E}_k$ for all $m \geq k$. 

\begin{lemm} \label{lemm:diameter_of_cube}
Let $x \in \partial Q_k(I)$ for some cube $Q_k(I) \in \mathcal{Q}_k$ and let $y \in Q_k^*(I)$. Then for all $m \geq k$, 
\begin{equation} 
d_m(x,y) \leq C_1r_k(I)|x-y|,
\end{equation} 
where 
$$C_1 = \frac{2nM(M-2n+1)R}{1-(M-2n+1)/M}.$$
\end{lemm}
\begin{proof}
First, we assume that $|x-y| \geq M^{-(k+1)}$. 

Let $x_k \in \mathcal{V}_k$ be a common vertex of $Q_k(I)$ and $Q_k(y)$, where we choose $Q_k(y)$ so that $Q_k(y) \subset Q_k^*(I)$. There exists a path $\widetilde{\gamma}$ in $\mathcal{E}_k \cap Q_k(I)$ connecting $x$ to $x_k$ with Euclidean length at most $nM^{-k}$. In fact, we can take $\widetilde{\gamma}$ to be a polygonal path whose segments are each parallel to a different coordinate axis. Observe that $\ell_m(\widetilde{\gamma}) \leq nr_k(I)M^{-k}$. 

Next, let $x_{k+1} \in \mathcal{V}_{k+1}$ be a point in $Q_{k+1}(y)$ (again, choose $Q_{k+1}(y)$ so that $Q_{k+1}(y) \subset Q_k(y)$). 
We similarly choose a polygonal path $\gamma_k$ which connects $x_k$ to $x_{k+1}$ and is contained in $\mathcal{E}_{k+1} \cap Q_k(y)$, with Euclidean length at most $nM^{-k}$. Observe that $\rho_m(w) \leq (M-2n+1)r_k(y)$ for all $w \in |\gamma_k|$. Hence $\gamma_k$ satisfies  $$\ell_m(\gamma_k) \leq n(M-2n+1)r_k(y)M^{-k}.$$ 	

We continue in a similar fashion. Let $x_{k+2} \in \mathcal{V}_{k+2}$ be a point in $Q_{k+2}(y) \subset Q_{k+1}(y)$. Define a polygonal path $\gamma_{k+1}$ from $x_{k+1}$ to $x_{k+2}$ in the same manner as the path $\gamma_k$. This path satisfies $$\ell_m(\gamma_{k+1}) \leq n(M-2n+1)r_{k+1}(y)M^{-(k+1)}.$$ 

Define points $x_j$ and paths $\gamma_j$ in the same manner until reaching $x_m$. Now, let $\gamma_m$ denote the straight-line path from $x_m$ to $y$, which satisfies $\ell_m(\gamma_m) \leq nr_m(x_m)M^{-m}$. The concatenation of $\widetilde{\gamma}, \gamma_k, \gamma_{k+1}, \ldots, \gamma_m$ is a path $\gamma$ from $x$ to $y$ satisfying 
\begin{align*} 
\ell_m(\gamma) & \leq nr_k(I)M^{-k} + n(M-2n+1) \sum_{j=k}^m r_j(y)M^{-j} \\
& \leq nr_k(I)M^{-k} + n(M-2n+1)r_k(y)M^{-k} \sum_{j=0}^{m-k} \left(\frac{M-2n+1}{M}\right)^j \\ 
& \leq nr_k(I)M^{-k} + n(M-2n+1)Rr_k(I)M^{-k} \left( \frac{1}{1-(M-2n+1)/M}\right).
\end{align*}
To obtain the middle inequality, we use the fact that $r_{j+1}(y) \leq (M-2n+1)r_j(y)$. The first term in the last line is smaller than the second term, so we obtain 
$$\ell_m(\gamma) \leq 2nR(M-2n+1)r_k(I)M^{-k} \left( \frac{1}{1-(M-2n+1)/M}\right).$$
Finally, since $M^{-k} \leq M|x-y|$, this establishes the result for the first case.

Next, suppose that $|x-y| < M^{-(k+1)}$. Choose $Q_{k+1}(x)$ so that $Q_{k+1}(x) \subset Q_k(I)$. Observe then that $y \in Q_{k+1}^*(x)$ and that $r_{k+1}(x) = r_k(I)$. Applying the first portion of the proof inductively yields the lemma for all $y \in Q_k^*(I)$. 
\end{proof}

The final lemma of this section deals with estimating $d_m(x,y)$ when the points $x,y \in Q$ lie at roughly distance $M^{-k}$ apart, for $m \geq k$. 

\begin{lemm} \label{lemm:estimate_for_k_separated_points}
Let $x,y \in Q$ be distinct, and let $k$ be the smallest value for which $Q_k(x) \cap Q_k(y) = \emptyset$. Then for all $m \geq k$, 
\begin{equation} 
\frac{r_k(x)|x-y|}{2nMR} \leq d_m(x,y) \leq C_1C_2r_k(x)|x-y|,
\end{equation}   	
where $C_1$ is the constant in Lemma \ref{lemm:diameter_of_cube} and
$$C_2 = \frac{R^{2nM+2}-1}{R-1}.$$ 
\end{lemm}
\begin{proof}
Observe that our hypotheses imply that $M^{-k} \leq |x-y| \leq 2nM^{-k+1}$. For the lower bound, we note that any path $\gamma$ from $x$ to $y$ must satisfy $$\ell_k(\gamma) \geq \frac{r_k(x)M^{-k}}{R} \geq \frac{r_k(x)|x-y|}{2nMR},$$ since it must cross $\overline{Q_k^*(x) \setminus Q_k(x)}$. By Lemma \ref{lemm:length_of_curves}, the same inequality holds with $\ell_m$ in place of $\ell_k$, for all $m \geq k$. The lower bound follows.  
	
Next, we obtain the upper bound. Let $x'$ be a vertex of $Q_k(x)$ and let $y'$ be a vertex of $Q_k(y)$. From the minimality of $k$, there is a polygonal path $\gamma$ from $x'$ to $y'$ formed by connecting points $(x_j)_{j=1}^{J}$ in $\mathcal{V}_k$, where $x_1 = x'$, $x_J = y'$, $d(x_{j-1}, x_j) = M^{-k}$ and $J \leq 2nM$. Also, let $x_0 = x$ and $x_{J+1} = y$. Using Lemma \ref{lemm:diameter_of_cube}, $d_m(x,x') + \ell_m(\gamma) + d_m(y',y)$ is bounded above by
\begin{align*} 
\sum_{j=0}^{J+1} C_1r_k(x_j)M^{-k} & \leq \sum_{j=0}^{J+1} C_1r_k(x)R^jM^{-k}  \\ 
& \leq C_1r_k(x) M^{-k}\frac{(R^{2nM+2}-1)}{R-1}.
\end{align*} 
Since $M^{-k} \leq |x-y|$, this gives the upper bound.  
\end{proof}

\section{Quasisymmetric equivalence} \label{sec:quasisymmetry}

We now formally define the metric space $(X,d)$ and the map $f\colon Q \rightarrow X$. Let $x,y \in \mathcal{V}$. Lemma \ref{lemm:increasing_property} implies that the sequence $(d_k(x,y))$ is eventually increasing. Lemma \ref{lemm:estimate_for_k_separated_points} implies that the sequence $(d_k(x,y))$ is bounded. Define a metric $d$ on $\mathcal{V}$ by $d(x,y) = \lim_{k \rightarrow \infty} d_k(x,y)$, which limit must exist for all $x,y \in \mathcal{V}$. As mentioned earlier, Lemma \ref{lemm:diameter_of_cube} implies that the metric $d$ is continuous with respect to the Euclidean distance on the set $\mathcal{V}$. 

We wish to extend $d$ continuously to all of $Q$. To verify that such an extension exists, we first show that for any two sequences $(x_j)$, $(x_j')$ in $\mathcal{V}$ converging to a point $x \in Q$, we have $d(x_j,x_j') \to 0$. Without loss of generality, we may assume that $x_j \neq x_j'$ for all $j \in \mathbb{N}$. For a given pair of points $x_j, x_j'$, let $k \in \mathbb{N}$ be as in the statement of Lemma \ref{lemm:estimate_for_k_separated_points}. Recall from the proof of Lemma \ref{lemm:estimate_for_k_separated_points} that $|x_j-x_j'| \leq 2nM^{-k+1}$. In addition, we have the universal upper bound $r_k(x_j) \leq (M-2n+1)^k$. Combining these estimates gives
$$d(x_j,x_j') \leq 2nMC_1C_2\left(\frac{M-2n+1}{M}\right)^k.$$
Since $M^{-k} \leq |x_j-x_j'|$, it follows that $k \to \infty$ as $j \to \infty$ and hence that the right-hand side of the above inequality converges to zero. 

We now define the metric $d$ on all of $Q$ by
$$d(x,y) = \lim_{j \to \infty} d(x_j,y_j),$$
where $(x_j)$ is a sequence in $\mathcal{V}$ converging to $x$ and $(y_j)$ is a sequence in $\mathcal{V}$ converging to $y$. It follows from the previous paragraph that the value of $d(x,y)$ is independent of the choice of sequence. We let $X$ denote the resulting metric space and take $f\colon Q \rightarrow X$ to be the identity map.


\begin{prop}
The map $f\colon Q \rightarrow X$ is quasisymmetric.
\end{prop} 
\begin{proof} 
Let $t>0$ be given. Consider an arbitrary triple of distinct points $x,y,z \in Q$ satisfying $|x-y| \leq t|x-z|$. As in Lemma \ref{lemm:estimate_for_k_separated_points}, let $k$ be the smallest value for which $Q_k(x) \cap Q_k(z) = \emptyset$. Then Lemma \ref{lemm:estimate_for_k_separated_points} gives $d(x,z) \geq r_k(x)|x-z|/(2nMR)$. As before, we observe that $M^{-k} \leq |x-z| \leq 2nM^{-k+1}$.  

Next, let $m$ be the smallest value for which $Q_m(x) \cap Q_m(y) = \emptyset$, noting that $M^{-m} \leq |x-y| \leq 2nM^{-m+1}$. Applying Lemma \ref{lemm:estimate_for_k_separated_points} gives $d(x,y) \leq C_1C_2 r_m(x)|x-y|.$ 
We split now into two cases. First, assume that $m \geq k$. Then $r_m(x) \leq (M-2n+1)^{m-k}r_k(x)$ and
\begin{align} \label{equ:m_geq_k_estimate}
d(x,y) & \leq 2nMC_1C_2 r_k(x)M^{-k}\left(\frac{M-2n+1}{M}\right)^{m-k} .
\end{align}
Since $M^{-k} \leq |x-z|$, we have
\begin{equation} \label{equ:quasisymmetric_bound} 
\frac{d(x,y)}{d(x,z)} \leq 4n^2M^2C_1C_2R \left( \frac{M-2n+1}{M}\right)^{m-k}. 
\end{equation}
The inequality $M^{-m} \leq |x-y| \leq t|x-z| \leq t2nM^{-k+1}$ now yields
$$m-k \geq  \max\left\{0, \frac{-\log(2nMt)}{\log M}\right\} .$$ 
Hence \eqref{equ:quasisymmetric_bound} gives a bound $\eta(t)$ on $d(x,y)/d(x,z)$ with the property that $\eta(t) \rightarrow 0$ as $t \rightarrow 0$.

In the case that $k \geq m$, instead of the estimate \eqref{equ:m_geq_k_estimate}, we have
\begin{equation*} \label{equ:m_leq_k_estimate}  
d(x,y) \leq 2nMC_1C_2M^{-m}r_k(x) M^{k-m} \leq 2nMC_1C_2r_k(x) M^{2(k-m)},
\end{equation*}
using the fact that $r_m(x) \leq M^{k-m}r_k(x)$. As before, this gives
\begin{equation} \label{equ:quasisymmetric_bound2} 
\frac{d(x,y)}{d(x,z)} \leq 4n^2M^2C_1C_2R M^{2(k-m)}. 
\end{equation}
Observe that $1 \leq 2nMt$, which gives
$$k-m \leq \frac{\log(2nMt)}{\log M}.$$ 
From here, \eqref{equ:quasisymmetric_bound2} now gives a bound $\eta(t)$ on $d(x,y)/d(x,z)$. 
\end{proof}

\section{Failure of inverse absolute continuity} \label{sec:absolute_continuity} 

The idea of this section is to identify each weight $\rho_k$ with a random variable $Y_k$ on $Q$ in the natural way, taking Lebesgue $n$-measure as the corresponding probability measure. Letting $X_k = \rho_k/\rho_{k-1}$ for $k \geq 1$, we obtain a sequence of independent, identically distributed random variables $(X_k)_k$ whose distribution is given by 
\begin{align*}
\mathbb{P}(X_k = 1) = & \frac{M^n - (M-2)^n}{M^n}, \\
\mathbb{P}(X_k = M-2n+1) = & \frac{(M-2)^n - (M-2n)^n}{M^n}, \\
\mathbb{P}(X_k = 1/L) = & \frac{(M-2n)^n}{M^n} .
\end{align*}
Each of the random variables $X_k$ has geometric mean 
$$\mu = (M-2n+1)^{((M-2)^n - (M-2n)^n)/M^n}\left(\frac{1}{L} \right)^{(M-2n)^n/M^n}, $$
where we consider $\mu$ as a function of $M$ and $L$. Since $Y_k = X_1X_2 \cdots X_k$, the strong law of large numbers implies that the sequence $(Y_k^{1/k})$ converges almost surely to $\mu$. Let $E_0 \subset Q$ be the set of points $x \in Q$ for which $Y_k^{1/k}(x) \rightarrow \mu$. We later define $E$ to be a certain full-measure subset of $E_0$. 

For the basic conclusion of Theorem \ref{thm:main} to hold (that is, without any quantitative statement about the Hausdorff dimension of the target), we need only choose $M$ and $L$ sufficiently large so that $\mu <1$. However, we will be more precise here. For the remainder of this section, we will fix a real number $\beta>0$ and let $L = M^\beta$. That is, we consider the value $L$ in our construction as a function of $M$.
The next lemma will allow us to estimate the Hausdorff dimension of $f(E)$.
\begin{lemm} \label{lemm:limit}
For all $\alpha \in (n/(1+ \beta), n)$, 
\begin{equation} \label{equ:limit}
\lim_{M \rightarrow \infty} M^{n-\alpha}\mu^\alpha = 0. 
\end{equation}
\end{lemm}
\begin{proof}
We observe first that
$$\lim_{M \rightarrow \infty} (M-2n+1)^{((M-2)^n - (M-2n)^n)/M^n}=1.$$
It follows from this that
\begin{align*} 
\lim_{M \rightarrow \infty} M^{n-\alpha}\mu^\alpha & = \lim_{M \rightarrow \infty} M^{n-\alpha} \left(\frac{1}{M} \right)^{\alpha\beta(M-2n)^n/M^n} \\
 & = \lim_{M \rightarrow \infty} M^{n-\alpha(1+ \beta)},
\end{align*}
where the limit on the left-hand side exists whenever the limit on the right-hand side exists. The latter limit is defined and equal to zero if $n-\alpha(1+ \beta)<0$. This holds if $\alpha \in (n/(1+ \beta), n)$.
\end{proof}

Now, we fix a value $\alpha \in (n/(1+ \beta), n)$. Observe that by initially choosing $\beta$ to be sufficiently large, we can choose $\alpha$ to be arbitrarily close to zero. The value $\alpha$ will be our bound on the Hausdorff dimension of $f(E)$.  Next, applying Lemma \ref{lemm:limit}, choose $M$ to be sufficiently large so that $M^{n-\alpha} \mu^\alpha < 1/2^\alpha$. 

Observe that $Y_k$ is constant on the interior of each cube $Q_k(I)$ of level $k$ with value $r_k(I)$, and that there are $M^{nk}$ such cubes. Let $m \in \mathbb{N}$. By the almost sure convergence of $(Y_k^{1/k})_k$ to $\mu$, we can pick $k_m$ sufficiently large so that $r_{k_m}(I) \leq (1+2^{-m})^{k_m} \mu^{k_m}$ holds for at least $(1-2^{-m})M^{nk_m}$ cubes $Q_{k_m}(I) \in \mathcal{Q}_{k_m}$ of level $k_m$ (out of $M^{nk_m}$ total cubes in $\mathcal{Q}_{k_m}$). We can also require that $k_m \geq m$. Let $\mathcal{F}_m \subset \mathcal{Q}_{k_m}$ denote the subcollection for which this holds, and let $F_m = \bigcup \mathcal{F}_m$. Finally, let 
$$E = \bigcap_{m=1}^\infty \bigcup_{j=m}^\infty F_j.$$
One checks that $[0,1]^n \setminus E$ has measure zero, and hence that $E$ is a set of full measure in $[0,1]^n$. 

Let $\varepsilon = \varepsilon(m) = 2nC_1(2\mu/M)^{k_m}$, where $C_1$ is the constant in Lemma \ref{lemm:diameter_of_cube}. We claim that for all $Q \in \mathcal{F}_m$, the set $f(Q)$ has $d$-diameter at most $\varepsilon$. By the triangle inequality applied to each coordinate function, we see that $|x-y| \leq nM^{-k_m}$ for all $x,y \in Q$. If it also holds that $y \in \partial Q$, Lemma \ref{lemm:diameter_of_cube} yields
\begin{align*} 
d(x,y) & \leq C_1r_{k_m}(I)nM^{-k_m} \\
 &  \leq nC_1 (1+2^{-m})^{k_m}\mu^{k_m} M^{-k_m} \\ 
 & \leq nC_1(2\mu/M)^{k_m}.
\end{align*}  
By fixing a base point on $\partial Q$ and applying the triangle inequality, it follows that $d(x,y) \leq 2nC_1(2\mu/M)^{k_m} = \varepsilon$ for all $x,y \in Q$. Thus $\diam f(Q) \leq \varepsilon$ as claimed. 

From this, we obtain the bound
\begin{align*} 
\mathcal{H}_\varepsilon^\alpha(F_m) & \leq \sum_{Q \in \mathcal{F}_m} (\diam f(Q))^\alpha \\ & \leq M^{nk_m} (2nC_1)^\alpha \left(2\mu/M \right)^{\alpha k_m} \\ 
  & \leq (2nC_1)^\alpha \left( M^{n-\alpha} (2\mu)^\alpha \right)^m,
\end{align*}
where in the final line we use the relationships $k_m \geq m$ and $M^{n-\alpha} (2\mu)^\alpha < 1$. 
Here, $\mathcal{H}_\varepsilon^\alpha$ denotes the $\alpha$-dimensional Hausdorff $\varepsilon$-content relative to the metric $d$. This shows in particular that 
\begin{align*} 
\mathcal{H}_\varepsilon^\alpha\left(E\right) & \leq (2nC_1)^\alpha \sum_{j=m}^\infty \left( M^{n-\alpha} (2\mu)^\alpha \right)^m  = (2nC_1)^\alpha \frac{\left( M^{n-\alpha} (2\mu)^\alpha \right)^m}{1- M^{n-\alpha} (2\mu)^\alpha}. 
\end{align*}
Since $M^{n-\alpha} (2\mu)^\alpha<1$, the expression on the right goes to zero as $m \rightarrow \infty$. We also observe that $\varepsilon(m) \to 0$ as $m \to \infty$; this suffices to show that $\mathcal{H}^\alpha(f(E)) = 0$. We conclude that $\alpha$ is an upper bound on the Hausdorff dimension of $f(E)$.

\section{Proof of Theorem \ref{thm:main2}} \label{sec:theorem2_proof} 

Theorem 1.2 can be obtained by modifying the construction of $X$ in Section \ref{sec:construction} in one detail. The idea is to place a cap on the size of the weights $\rho_k$. For this section we assume that the parameter $L$ in Section \ref{sec:construction} is chosen to satisfy $L \geq M$.

As in Section \ref{sec:absolute_continuity}, we will interpret our construction from a probabilistic viewpoint. However, we do so in a somewhat different way, with the goal of keeping the analysis as simple as possible. For all $k \in \mathbb{Z}_{\geq 0}$, let $X_k$ be the random variable on $Q$ defined by
$$X_k(x) = \left\{ \begin{array}{cc} 1 & \text{ if } x \in P_k^1(Q_k(x)) \cup P_k^2(Q_k(x)) \\ -1 & \text{ if } x \in P_k^3(Q_k(x)) \end{array} \right. . $$
The sets $P_k^j(Q_k(x))$ ($j \in \{1,2,3\}$) are defined in the same way as the sets $P_k^j(I)$ in Section \ref{sec:construction}. Here, we require that the sets $Q_k(x)$ are chosen so that $Q_{k+1}(x) \subset Q_k(x)$ for all $k \in \mathbb{Z}_{k \geq 0}$, so that $X_k$ is well-defined for all $x \in Q$.  Observe that the random variables $X_k$ are independent and identically distributed. Next, for all $k \in \mathbb{Z}_{\geq 0}$, let $Y_k = X_0 + \cdots + X_k$. In this manner, we obtain for each $x \in Q$ a corresponding random walk $(Y_k)_k$ on the integers beginning at zero. This random walk steps up with probability $p=1 - (M-2n)^n/M^n$ and steps down with probability $q=(M-2n)^n/M^n$.

Now we explain the modification. In the course of applying the inductive definition \eqref{equ:rho_definition}, if there exits $k$ such that $Y_k(x) = 1$, then modify the inductive definition \eqref{equ:rho_definition} to set $$\rho_m(x) = \rho_k(q_k(I))$$ for all $m \geq k$ and interior points $x$ of $Q_k(I)$. Extend the definition of $\rho_m$ to all of $Q$ by lower semicontinuity. Let $(X,d)$ denote the resulting limit metric space, defined as before except with this modification. The proof that $(X,d)$ is quasisymmetrically equivalent to $Q$ still applies without modification. 

It is a basic exercise in probability to show that, for $q>1/2$, the probability $r$ that $Y_k(x)=1$ for some $k \in \mathbb{Z}_{\geq 0}$ is $r=(1-q)/q$. To derive this, one first notes that $r$ satisfies the recursive relationship $r = p + qr^2$. Solving for $r$ yields the solutions $r=p/q = (1-q)/q$ and $r=1$. However, the possibility that $r=1$ is ruled out by the transience of biased random walks, since $Y_k/k$ converges almost surely to $p-q$. See, for example, the book of Klenke \cite[Sec. 17.5]{Kle:08}.

Let $F$ be the set of points $x \in Q$ for which $Y_k(x) \leq 0$ for all $k \in \mathbb{Z}_{\geq 0}$. We have shown that $|F| = 1-(1-q)/q = (2q-1)/q$. Observe that as $M \rightarrow \infty$, $q \rightarrow 1$, so that $|E|$ can be made arbitrarily close to 1. Let $E$ be as in Section \ref{sec:absolute_continuity}. Then $E' = E \cap F$ is a set of full measure in $F$. By choosing $M$ and $L$ sufficiently large, as in Section \ref{sec:absolute_continuity}, we may ensure that $f(E')$ has Hausdorff dimension arbitrarily close to 0. 

To conclude the proof, we note that $f$ is a Lipschitz function with Lipschitz constant $(M-2n+1)$, since we have assumed that $L \geq M$. This implies that $X$ has finite Hausdorff $n$-measure. \\

\noindent {\bf Acknowledgments.} I am grateful to Kai Rajala, Tuomas Orponen, David Bate, and Enrico Le Donne for conversations about the content of this paper. Specifically, I thank Tuomas Orponen and David Bate for their observation that $f(E)$ in the construction has Hausdorff dimension strictly smaller than $n$. I also thank Kai Rajala and Tuomas Orponen for feedback on a draft of this paper, and the referee for additional comments and corrections. 

\
\bibliographystyle{abbrv}
\bibliography{biblio} 
\end{document}